\documentclass[12pt]{article}
\usepackage{bbm}
\usepackage{amsfonts}
\usepackage{pifont}
\usepackage{hyperref}
\usepackage{mathrsfs}
\usepackage{indentfirst}
\usepackage{amsmath}
\usepackage{amssymb}
\usepackage[amsmath, thmmarks]{ntheorem}
\usepackage{cite}
\usepackage{graphicx}
\usepackage{tikz}
\newtheorem{theorem}{Theorem}[section]
\newtheorem{definition}[theorem]{Definition}

\newtheorem{remark}[theorem]{Remark}
\newtheorem{corollary}[theorem]{Corollary}
\newtheorem{proposition}[theorem]{Proposition}
\newtheorem{example}[theorem]{Example}

\def\QEDopen{{\setlength{\fboxsep}{0pt}\setlength{\fboxrule}{0.2pt}\fbox{\rule[0pt]{0pt}{1.3ex}\rule[0pt]{1.3ex}{0pt}}}} 
\def\QED{\QEDopen}
\def\proof{{\bf Proof.} }
\def\endproof{\hspace*{\fill}~\QED\par\endtrivlist\unskip}

\begin{document}
\setcounter{page}{1}

\title{{\textbf{The pseudocomplementedness of modular lattices described by two 0-sublattices}}\thanks {The research of the first author was supported by the National Natural Science Foundation
of China No. 11901064. The research of the second author was supported by the National
Natural Science Foundation of China No. 12071325}}
\author{Peng He$^1$\footnote{\emph{E-mail address}: hepeng@cuit.edu.cn}, Xue-ping Wang$^2$\footnote{Corresponding author. xpwang1@hotmail.com; fax: +86-28-84761502}\\
\emph{\small 1. College of Applied Mathematics, Chengdu University of Information Technology}\\ \emph{\small Chengdu 610225, Sichuan, People's Republic of China}\\
\emph{\small 2. School of Mathematical Sciences, Sichuan Normal University}\\ \emph{\small Chengdu 610066, Sichuan, People's Republic of China}}

\newcommand{\pp}[2]{\frac{\partial #1}{\partial #2}}
\date{}
\maketitle

\begin{quote}
{\bf Abstract} In this article, we first characterize pseudocomplemented inductive modular lattices by using their two 0-sublattices. Then we use two 0-sublattices of a subgroup lattice to describe all locally cyclic abelian groups. In particular, we show that a locally cyclic abelian group can be characterized by its three subgroups.\\


\end{quote}

\section{Introduction}\label{intro}
In this article, we are interested in the lattice formed by the totality of subgroups
(resp. normal subgroups) of a group. Defining the meet , $\wedge$, and the join, $\vee$, of
subgroups of a group $G$
in the natural way, one sees easily that the totality of subgroups (resp. normal subgroups) of a
group $G$ satisfies all axioms of a lattice. We shall call this lattice a
subgroup (resp. normal subgroup) lattice of $G$ and denote it by $L(G)$ (resp. $N(G)$).
We first recall the well-known facts of lattices $L(G)$ and $N(G)$ as follow.
\begin{theorem}[\cite{Birkhoff73,Gratzer14,Suz}]\label{theo001}
For any group $G$, $L(G)$ is an algebraic (compactly generated) lattice, and $N(G)$ is a modular lattice
and  $A\vee B=AB=BA$ for all normal subgroups $A$ and $B$ of $G$ .
\end{theorem}

In studying various algebraic systems emphasis has been put on
the structure of their subsystems rather than on the behavior of individual elements
in the systems, and the lattice of subgroups has drawn
the particular attention of mathematicians since the birth of group theory.
It was Dedekind who considered the system of ideals in a ring
of algebraic integers for the first time from the lattice theoretical point
of view, and he discovered and used the modular identity, sometimes
called the Dedekind law, in his calculation of ideals. But the real
history of the theory of subgroup lattices began in 1928, when Rottl\"{a}nder
considered in her paper \cite{Rott} the totality of subgroups of
a finite group and the mappings between subgroup lattices in solving a
question arising from field extensions. We refer to Suzuki's book\cite{Suz}
and Schmidt's book\cite{Schm}
for more information about those theory.

Let $L=(L, \vee, \wedge)$ be a lattice with the least element $0$ and the greatest element $1$, and let $a\in L$.
An element $b\in L$ is called a complement of $a$ if $a\vee b=1$ and $a\wedge b=0$, and the lattice
$L$ is complemented if every element of $L$ has a complement. An element $a^*\in L$ is called
a pseudocomplement of $a$ if  $a\wedge a^*=0$ and $a\wedge x=0$ implies $x\leq a^*$.
The element $a^*$ is obviously the greatest element in the set of all $x\in L$ for which $a\wedge x=0$.
In other words, the subset of all elements disjoint from $x$ is required to form a principal ideal.
A lattice is called pseudocomplemented if each of its elements has a pseudocomplement. Note
that the terminology is slightly misleading, since a complement is not necessarily a
pseudocomplement (see, e.g., \cite{Stern99, Crawley73, Frink}).

It is an interesting question in group theory to what extent the structure of
the subgroup lattice of a group determines the structure of the group itself.
Suzuki spent his early research years on this problem \cite{Suz51,Suz52}.
Since then, many characterizations and classifications
have been obtained for groups for which the subgroup lattice or normal subgroup
lattice has certain lattice-theoretic properties. Possibly the most famous result in this
direction is Ore's result that a group is locally cyclic if and only if its lattice of
subgroups is distributive \cite{Ore}. C\u{a}lug\u{a}reanu \cite{Calu86} further proved the following theorem if the group is abelian:
\begin{theorem}[\cite{Calu86}]\label{theo002}
For an abelian group $G$ the following three conditions are equivalent: \emph{(i)} $L(G)$ is a distributive lattice; \emph{(ii)}  $L(G)$ is pseudocomplemented; \emph{(iii)} $G$ is a locally cyclic group.
\end{theorem}
 He also characterized abelian groups which have a Stone lattice or a Heyting algebra of subgroups (see \cite{Gratzer14, Chen} for the detailed investigations of Stone lattices and Heyting algebras). Recently, Medts and T\u{a}rn\u{a}uceanu \cite{Metd, Tarn04, Tarn05} studied finite groups by admitting a pseudocomplemented
subgroup lattice or a pseudocomplemented normal subgroup lattice.

It is worth to be pointed out that Katri\v{n}\'{a}k and Mederly \cite{Katrinak74}
investigated modular pseudocomplemented lattices in terms of triples, and Chameni-Nembua
and Monjiardet \cite{Chameni} gave a simple criterion for pseudocomplementedness of strongly atomic
algebraic lattices by its atoms (also see \cite{Stern99}). In particular, He and Wang \cite{He21} characterized pseudocomplemented lattices
with (ACC) and (DCC) by using their nine 0-sublattices that are defined as follows.
\begin{definition}[\cite{He21}]\label{de1}
\emph{We say that a sublattice of a lattice $L$ with the least element is a $0$-sublattice if it contains the least element of $L$.}
\end{definition}
By using Dedekind's modularity criterion that a lattice is modular if and only if it contains no
five-element sublattice isomorphic to $N_5$,
they proved that a modular lattice $L$ with (ACC) and (DCC) is pseudocomplemented  if and only if $L$
contains no 0-sublattice isomorphic with one of the lattices $M_3$ and $M_{2,3}$
represented in Figure 1.
\par\noindent\vskip10pt
\begin{minipage}{11pc}
\setlength{\unitlength}{0.75pt}\begin{picture}(600,160)
\put(180,140){\circle{6}}
\put(180,60){\circle{6}}
\put(140,100){\circle{6}}
\put(180,100){\circle{6}}
\put(220,100){\circle{6}}
\put(180,63){\line(0,1){34}}
\put(180,103){\line(0,1){34}}
\put(178,62){\line(-1,1){36}}
\put(182,62){\line(1,1){36}}
\put(142,102){\line(1,1){36}}
\put(218,102){\line(-1,1){36}}
\put(170,42){$M_3$}

\put(380,40){\circle{6}}
\put(380,160){\circle{6}}
\put(380,80){\circle{6}}
\put(340,120){\circle{6}}
\put(380,120){\circle{6}}
\put(420,120){\circle{6}}
\put(420,80){\circle{6}}
\put(380,83){\line(0,1){34}}
\put(420,83){\line(0,1){34}}
\put(380,43){\line(0,1){34}}
\put(380,123){\line(0,1){34}}
\put(378,82){\line(-1,1){36}}
\put(382,82){\line(1,1){36}}
\put(382,42){\line(1,1){36}}
\put(342,122){\line(1,1){36}}
\put(418,122){\line(-1,1){36}}
\put(370,22){$M_{2,3}$}
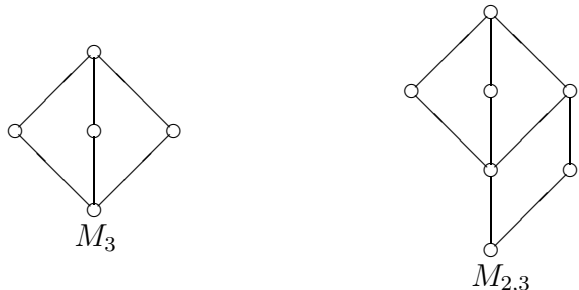
\put(115,3){FIGURE 1. Two nonpseudocomplemented lattices}
\end{picture}
\end{minipage}

As for each finite abelian group $G$, the subgroup lattice $L(G)$ equals to the normal subgroup lattice $N(G)$. Combining Theorems \ref{theo001} and \ref{theo002} with the result of He and Wang in \cite{He21} mentioned as above, we get
the following corollary obviously.

\begin{corollary}\label{coro1}
For a finite abelian group $G$, the following four conditions are equivalent:\\ \emph{(i)} $L(G)$ is a distributive lattice;\\ \emph{(ii)}  $L(G)$ is pseudocomplemented;\\ \emph{(iii)} $G$ is a locally cyclic group;\\
\emph{(iv)} $L(G)$ contains no 0-sublattice isomorphic to one of the lattices $M_3$ and $M_{2,3}$.
\end{corollary}

Motivated by the characterizations of pseudocomplemented modular lattice with (ACC) and (DCC),
in this article, we first characterize the pseudocomplementedness of inductive modular lattices
by using their 0-sublattices equivalently. Then we use two 0-sublattices of a subgroup lattice to describe all locally cyclic abelian groups.

\section{Main results}
We assume that the reader is familiar with the basic
theory of lattices such as a partially ordered set (poset), a chain, a
lattice, a subgroup lattice, a normal subgroup lattice, etc.,
(see, e.g., \cite{Birkhoff73,Crawley73, Suz, Schm}).

\begin{definition}[\cite{Calu}]\label{de2}
\emph{We say that a lattice $L$ is inductive if each quotient sublattice
$[a,b]$ of it satisfies the following condition (B).\\
(B) if for any chain $\{b_i\}_{i\in I}$ in $[a,b]$ and for any $x\in [a,b]$ such that
$x\wedge b_i=a$ for all $i\in I$, there holds $x\wedge \bigvee_{i\in I}b_i=a$.}
\end{definition}

Obviously, intervals of inductive lattices are inductive and each upper continuous lattice is inductive.
Thus, each algebraic lattice is also inductive since it is upper continuous (see \cite{Crawley73,Birkhoff73}).

Let $L$ be a lattice with the least element $0$ and the greatest element $1$ (No problem arises in distinguishing these elements from the integers $0$ and $1$). For all $a, b\in L$, $a\parallel b$ denotes that $a\ngeq b$ and
$a\nleq b$, and $a\nparallel b$ denotes that $a\geq b$ or $a\leq b$.
 Before formulating
our main results, we need two propositions as below.

\begin{proposition}\label{pro1}
Let $L$ be a modular lattice. If $b$ is a maximal element with $a\wedge b=0$ in $L$, then $(a\vee b)\wedge x\neq 0$ for all $x\neq 0$.
\end{proposition}
\proof
Let $c>b$ in $L$. We claim that $b<b\vee (a\wedge c)$.
Otherwise, $a\wedge c\leq b$, this implies that
$a\wedge c=(a\wedge c) \wedge b=c\wedge (a\wedge b)=0$.
However, $a\wedge c\neq 0$ since $b$ is a maximal element with $a\wedge b=0$ and $c>b$, a contradiction.
As $L$ is modular and $b<c$, $b\vee (a\wedge c)=(a\vee b)\wedge c$. Thus
\begin{equation}\label{eq1}
 b<b\vee (a\wedge c)=(a\vee b)\wedge c.
\end{equation}

Suppose that $x\neq 0$. If $x\leq b$, then $(a\vee b)\wedge x=x\neq 0$. If $x\nleq b$, then $b<b\vee x$. Thus  $b<(a\vee b)\wedge (b\vee x)=((a\vee b)\wedge x)\vee b$ by using
formula (\ref{eq1}) and the modularity of $L$. Therefore, $(a\vee b)\wedge x\neq 0$ for all $x\neq 0$.
\endproof

\begin{proposition}\label{pro2}
Let $L$ be an inductive lattice. Then there is at least one maximal element $b$ such that $a\wedge b=0$ in $L$.
\end{proposition}
\proof
Set $\mathcal{S}=\{x\in L\mid a\wedge x=0\}$. Obviously, $\mathcal{S}$ is a nonempty
set of $L$ since $0\in \mathcal{S}$. $L$ being inductive, one can easily check that
each chain $C$ in $\mathcal{S}$ contains the element $\bigvee C$. Hence, by Zorn's Lemma, $\mathcal{S}$
has at least one maximal element $b$ satisfying $a\wedge b=0$.
\endproof

\begin{theorem}\label{theo1}
In an inductive modular lattice $L$, the following conditions are equivalent:\\
\emph{(a)} $L$ is pseudocomplemented;\\
\emph{(b)} $L$ contains no 0-sublattice isomorphic to one of the lattices $M_3$ and $M_{2,3}$;\\
\emph{(c)} $L$ contains no ternary sequence $(a, b, c)$ satisfying the following three conditions:

\emph{(i)} $0\notin \{a,b,c\}$; \emph{(ii)} $c\wedge a=c\wedge b=0$; \emph{(iii)} $c\vee a=c\vee b=a\vee b$.
\end{theorem}
\proof
(a) $\Rightarrow$ (b) $L$ being pseudocomplemented, the subset of all elements $x$ with $x\wedge a=0$ is required to form an ideal of $L$ for any element $a\in L$. Thus one can check that $L$ contains no 0-sublattice isomorphic to one of the lattices $M_3$ and $M_{2,3}$.

(b) $\Rightarrow$ (c) Suppose that there exists a ternary sequence $(a, b, c)$ satisfying conditions (i), (ii) and (iii). We claim that $a \parallel b$, $a \parallel c$ and $b \parallel c$. Obviously $a \parallel c$ and $b \parallel c$ since $0\notin \{a, b, c\}$ and $c\wedge a=c\wedge b=0$. If $a\nparallel b$, say $a<b$, then $a\vee c=b\vee c=a\vee b=b$. Thus $c\leq b$ and $0=b\wedge c=c$, contrary to the fact $c\neq 0$. The proof is made in two cases.

Case 1. If $a\wedge b=0$, then the elements $0, a, b, c$ and $a\vee c$ form a 0-sublattice of $L$ which is isomorphic to $M_3$.

Case 2. If $a\wedge b\neq 0$, then $0< a\wedge b<a<a\vee b$ since $a\parallel b$. If $c\vee (a\wedge b)=a\vee b$, then the elements $0, a\wedge b, a, a\vee b$ and $c$ form a nonmodular five-element sublattice $N_5$ of $L$, a contradiction. On the other hand, we know that $c< c\vee (a\wedge b)$ since $c\wedge (a\wedge b)=0$, $a\wedge b\neq 0$ and $c\neq 0$. Thus $c<c\vee (a\wedge b)<a\vee b$. In the meantime, by using modularity, we have
\begin{equation*}
[c\vee (a\wedge b)]\wedge a=(a\wedge b)\vee (a\wedge c)=a\wedge b \ \,(\mbox{since } a\wedge c=0)
\end{equation*}
and
\begin{equation*}
[c\vee (a\wedge b)]\wedge b=(a\wedge b)\vee (b\wedge c)=a\wedge b \ \,(\mbox{since } b\wedge c=0).
\end{equation*}
Therefore, the elements $0, a, b, c, a\wedge b, c\vee (a\wedge b)$ and $a\vee c$ form a 0-sublattice of $L$ which is isomorphic to $M_{2,3}$.

From Cases 1 and 2, the proof of (b) $\Rightarrow$ (c) is finished.

(c) $\Rightarrow$ (a) Now, suppose that $L$ is not a pseudocomplemented inductive modular lattice.
Then there exists an element $c$ ($c\neq 0$) that does not have a pseudocomplement. Thus, from
Proposition \ref{pro2}, there exist two different maximal elements $a$ and $b$ in $L$ with $c\wedge a=0$ and
$c\wedge b=0$, respectively. Obviously, $a\parallel b$, this yields that $a\neq 0$ and $b\neq 0$.
Since $a$ and $b$ play a symmetric role, we are going to list these cases only modulo $a-b$ symmetry.
It should be pointed out that the next proof involves a basic technique for constructing a
ternary sequence $(x, y, z)$ satisfying conditions (i), (ii) and (iii) in each case,
described in A, B and C as follows.

A. $a\vee c=b\vee c$. In fact, $a\vee c=b\vee c=a\vee b\vee c\geq a\vee b$.
Thus the proof is made in two steps.

A1. $a\vee c=b\vee c=a\vee b$. One can easily check that $(a, b, c)$ satisfies conditions (i), (ii) and (iii).

A2. $a\vee c=b\vee c>a\vee b$. From the maximal property of $a$ and $b$, $c\wedge (a\vee b)\neq 0$.
Using modularity of $L$, we have
\begin{equation*}
[c\wedge (a\vee b)]\vee a=(a\vee b)\wedge (a\vee c)=a\vee b \ \, (\mbox{since } a\vee b<a\vee c)
\end{equation*}
and
\begin{equation*}
[c\wedge (a\vee b)]\vee b=(a\vee b)\wedge (b\vee c)=a\vee b \ \, (\mbox{since } a\vee b<b\vee c).
\end{equation*}
Moreover, $[c\wedge (a\vee b)]\wedge a=c\wedge a=0$ and $[c\wedge (a\vee b)]\wedge b=c\wedge b=0$.
Consequently, the ternary sequence $(a, b, c\wedge(a\vee b))$ satisfies conditions (i), (ii) and (iii).

B. $a\vee c< b\vee c$. Then $a\vee b\leq b\vee c$ since $b\vee c=a\vee b\vee c$. Thus there are two subcases.

B1. $a\vee b=b\vee c$. By using Proposition \ref{pro1}, $(a\vee c)\wedge b\neq 0$.
Using modularity of $L$, we have
\begin{equation*}
[(a\vee c)\wedge b]\vee a=(a\vee c)\wedge (a\vee b)=a\vee c \ \,(\mbox{since } a\vee c<a\vee b)
\end{equation*}
and
\begin{equation*}
[(a\vee c)\wedge b]\vee c=(a\vee c)\wedge (b\vee c)=a\vee c \ \,(\mbox{since } a\vee c<b\vee c).
\end{equation*}
In the meantime, $a\wedge c=0$ and $[(a\vee c)\wedge b]\wedge c=c\wedge b=0$. Thus the ternary
sequence $(a, (a\vee c)\wedge b, c)$ satisfies conditions (i), (ii) and (iii).

B2. $a\vee b<b\vee c$. By using Proposition \ref{pro1} again, $(a\vee c)\wedge b\neq  0$.
Similarly to the proof of A2, we know that $(a\vee b)\wedge c\neq 0$. In the meantime, it is clear that
$[(a\vee b)\wedge c]\wedge a=a\wedge c=0$ and $[(a\vee b)\wedge c]\wedge [(a\vee c)\wedge b]=b\wedge c=0$.
As $L$ is modular, we obtain the following three equations:
\begin{equation*}
[(a\vee b)\wedge c]\vee a=(a\vee b)\wedge (a\vee c),
\end{equation*}

\begin{equation*}
[(a\vee c)\wedge b]\vee a=(a\vee b)\wedge (a\vee c)
\end{equation*}
and
\begin{align*}
[(a\vee b)\wedge c]\vee[(a\vee c)\wedge b]&=(a\vee c)\wedge [[(a\vee b)\wedge c]\vee b]\\
&=(a\vee c)\wedge (a\vee b)\wedge (b\vee c)\\
&=(a\vee c)\wedge(a\vee b) \ \,(\mbox{since } a\vee c<b\vee c).
\end{align*}
Hence the ternary sequence $(a, (a\vee c)\wedge b, (a\vee b)\wedge c)$ satisfies conditions (i), (ii) and (iii).

C. $a\vee c\parallel b\vee c$. We divide the proof of C into two subcases.

C1. $(a\vee c)\wedge (b\vee c)\leq a\vee b$. Again, by using Proposition \ref{pro1}, we have $a\wedge (b\vee c)\neq 0$ and $b\wedge (a\vee c)\neq 0$. Notice that
$c\wedge [a\wedge (b\vee c)]=c\wedge a=0$ and $c\wedge [b\wedge (a\vee c)]=c\wedge b=0$.
Since $L$ is modular, we have the following three formulas:
\begin{equation*}
[a\wedge (b\vee c)]\vee c=(b\vee c)\wedge (a\vee c),
\end{equation*}

\begin{equation*}
[b\wedge (a\vee c)]\vee c=(a\vee c)\wedge (b\vee c)
\end{equation*}
and
\begin{align*}
[a\wedge (b\vee c)]\vee [b\wedge (a\vee c)]&=(a\vee c)\wedge [[a\wedge (b\vee c)]\vee b]\\
&=(a\vee c)\wedge (a\vee b)\wedge (b\vee c)\\
&=(a\vee c)\wedge(b\vee c) \ \,(\mbox{since } (a\vee c)\wedge(b\vee c)\leq a\vee b).
\end{align*}
Consequently, the ternary sequence $(a\wedge (b\vee c), b\wedge (a\vee c), c)$
satisfies conditions (i), (ii) and (iii).

C2. $(a\vee c)\wedge (b\vee c)\nleq a\vee b$. First, we can easily check that $(a\vee b)\wedge c\neq 0$, $[a\wedge (b\vee c)]\wedge [(a\vee b)\wedge c]=a\wedge c=0$ and
$[b\wedge (a\vee c)]\wedge [(a\vee b)\wedge c]=b\wedge c=0$. Then by using the modularity of $L$, we have
\begin{equation*}
 [a\wedge (b\vee c)]\vee [b\wedge (a\vee c)]=(a\vee c)\wedge [[a\wedge (b\vee c)]\vee b]
=(a\vee c)\wedge (b\vee c)\wedge(a\vee b),
\end{equation*}
\begin{equation*}
[a\wedge (b\vee c)]\vee [(a\vee b)\wedge c]=(a\vee b)\wedge [[a\wedge (b\vee c)]\vee c]
=(a\vee b)\wedge (b\vee c)\wedge (a\vee c)
\end{equation*}
and
\begin{equation*}
[b\wedge (a\vee c)]\vee [(a\vee b)\wedge c]=(a\vee b)\wedge [[b\wedge (a\vee c)]\vee c]
=(a\vee b)\wedge (a\vee c)\wedge (b\vee c).
\end{equation*}
Similarly to the proof of C1, we know that $a\wedge (b\vee c)\neq 0$ and $b\wedge (a\vee c)\neq 0$.
Therefore, the ternary sequence $(a\wedge (b\vee c), b\wedge (a\vee c), (a\vee b)\wedge c)$ satisfies
conditions (i), (ii) and (iii).

A, B and C combine to prove (c) $\Rightarrow$ (a).
\endproof

\begin{remark}\label{remark11}
\emph{In 1900, Dedekind, in his famous paper \cite{Dedekind}, gave a three elements generated
free modular lattice $M_{28}$ that is a 28-element lattice (the diagram of $M_{28}$ can also be found in \cite{Birkhoff73}). Applying this result, one can prove that any modular
lattice generated by three elements $a, b$ and $c$ satisfying the condition (c) of Theorem \ref{theo1} must be either $M_3$ or $M_{2,3}$. Thus the conditions (b) and (c) in Theorem \ref{theo1} are equivalent.}
\end{remark}

Combining Theorem \ref{theo001} with Theorem \ref{theo1} we get:

\begin{corollary}\label{the001}
Let $(G, \cdot, e)$ be a group. If $N(G)$ is
inductive then the following three conditions are equivalent:\\
\emph{(i)} $N(G)$ is pseudocomplemented;\\
\emph{(ii)} $N(G)$ contains no 0-sublattice isomorphic to one of the lattices $M_3$ and $M_{2,3}$;\\
\emph{(iii)} $G$ contains no three normal subgroups $U, V, W$ different from $(\{e\}, \cdot, e)$
with $$U\cap W=V\cap W=\{e\} \mbox{ and }U\vee V=U\vee W=V\vee W \ (i.e.\,\, UV=UW=VW).$$
\end{corollary}

Combining Theorems  \ref{theo001} and \ref{theo002} and Corollary \ref{the001} with the fact that an algebraic lattice is inductive and $L(G)=N(G)$ for every abelian group $G$, we have the following theorem.

\begin{theorem}\label{the002}
Let $(G,\cdot, e)$ be a abelian group. The following five conditions are equivalent:\\
\emph{(i)} $L(G)$ is a distributive lattice;\\
\emph{(ii)} $G$ is a locally cyclic group;\\
\emph{(iii)} $L(G)$ is pseudocomplemented;\\
\emph{(iv)} $L(G)$ contains no 0-sublattice isomorphic to one of the lattices $M_3$ and $M_{2,3}$;\\
\emph{(v)} $G$ contains no three subgroups $U, V, W$ different from $(\{e\}, \cdot, e)$
with $$U\cap W=V\cap W=\{e\} \mbox{ and }U\vee V=U\vee W=V\vee W \ (i.e.\,\, UV=UW=VW).$$
\end{theorem}

Next, we shall give an example to illustrate Theorem \ref{the002}, which also tells us if the assumption of
inductive is dropped then Theorem \ref{theo1} fails.

\begin{example}\label{exa1}
\emph{Let $G$ be a cyclic group. Then $L(G)$ is pseudocomplemented.
In Particular, if $G$ is a finite cyclic group of order $n$ (resp. an infinite
cyclic group), then the lattice $L(G)$ is isomorphic to the lattice $L_n$
of all divisors of $n$ (resp. the lattice $(\mathbb{N}, \leq)$,
where $\mathbb{N}=\{0,1,\cdots,n,\cdots\},\mbox{ and }n_1\leq n_2$ ($n_1, n_2\in \mathbb{N}$)
if and only if $n_2|n_1$).}

\quad\\
\emph{(1) Suppose that there are three elements $n_1, n_2, n_3\in L_n$ and each of them is not equal to $n$ satisfying $n_1\wedge n_3=n_2\wedge n_3=n$ and $n_1\vee n_3=n_2\vee n_3=n_1\vee n_2$. Then  $n=\mbox{lcm}(n_1, n_3)=\mbox{lcm}(n_2, n_3)$ and $\mbox{gcd}(n_1, n_3)=\mbox{gcd}(n_2, n_3)=\mbox{gcd}(n_1, n_2)$.
Thus $n_1=n_2=\mbox{gcd}(n_1, n_2)$, this means that $n_1| n_3$ and $n_2|n_3$. Therefore, $n_3=n$, a contradiction.}

\emph{Now, we shall discuss the lattice $(\mathbb{N}, \leq)$. Notice that  $m\wedge n=0$ ($m, n\in \mathbb{N}$) if and only if $m=0$ or $n=0$. Thus $\mathbb{N}$ contains no three elements $n_1, n_2, n_3$ different from $0$ with $n_1\wedge n_3=n_2\wedge n_3=0$ and $n_1\vee n_3=n_2\vee n_3=n_1\vee n_2$.}

\quad\\
\emph{(2) Let us consider the dual lattice $(\mathbb{N}, \geq)$. One can check that the lattice $(\mathbb{N}, \geq)$ is modular since $(\mathbb{N}, \leq)$ is modular but $(\mathbb{N}, \geq)$ is not an inductive lattice. If there are three elements $n_1, n_2, n_3\in \mathbb{N}$ with $1\notin \{n_1, n_2, n_3\}$, $n_1\wedge n_3=n_2\wedge n_3=1$ ($1$ is the least element of $(\mathbb{N}, \geq)$) and $n_1\vee n_3=n_2\vee n_3=n_1\vee n_2$, then $\mbox{gcd}(n_1, n_3)=\mbox{gcd}(n_2, n_3)=1$ and $n_1 n_3=n_2 n_3=\mbox{lcm}(n_1, n_3)=\mbox{lcm}(n_2, n_3)=\mbox{lcm}(n_1, n_2)$. Thus $n_1=n_2=\mbox{lcm}(n_1, n_2)=n_1 n_3=n_2 n_3$, which yields that
$n_3=1$, a contradiction. Therefore,  $(\mathbb{N}, \geq)$ contains no three elements $n_1, n_2, n_3$ different from $1$ with $n_1\wedge n_3=n_2\wedge n_3=1$ and $n_1\vee n_3=n_2\vee n_3=n_1\vee n_2$.
However, there exist two elements $p, q\in \mathbb{N}$ such that $p\wedge q^k=\mbox{gcd}(p, q^k)=1$ for all $k\in \mathbb{N}$ (for example $p=3$ and $q=2$). This means that
there is no largest element $m$ in $(\mathbb{N}, \geq)$ such that $p\wedge m=1$. Therefore, the assumption of
inductive of the lattice $L$ in Theorem \ref{theo1} cannot be removed.}

\end{example}

\begin{remark}\label{re01}
\emph{(1) For an algebraic system $\Lambda$, the lattice $L(\Lambda)$ formed by the totality of subalgebras of $\Lambda$ is an algebraic lattice. Therefore, if $L(\Lambda)$ is modular, then $L(\Lambda)$ is pseudocomplemented if and only if $\Lambda$ contains no three subalgebras $A, B$ and $C$ different from the least subalgebra $O$ with $A\cap B=A\cap C=O \mbox{ and }A\vee B=A\vee C=B\vee C$.}

\emph{(2) Any sublattice of a modular lattice is modular. Then both the subspaces of any vector space $V$ and the
ideals of any ring $R$ form modular lattices being sublattices of the modular lattices of all (normal) subgroups
of the relevant additive group. Therefore, with a ring $R$ instead of a group $G$ and an ideal instead of a normal subgroup, and with a vector space $V$ instead of an abelian group $G$ and a subspace instead of a subgroup, respectively, Corollary \ref{the001} and Theorem \ref{the002} also hold.}
\end{remark}
\section{Conclusions}
This article gave us the characterizations of pseudocomplemented inductive modular lattices by using two 0-sublattices, which were applied to describe all locally cyclic abelian groups. As it is well known, if a group $G$ is not an abelian group then its subgroup lattice $L(G)$ is not a modular lattice generally. So that a challenge problem is whether we can move the modularity of the lattice.
\section*{Statements and Declarations}
\quad\\
{\bf{Ethical approval}} This article does not contain any studies with human
participants or animals performed by any of the authors.\\
{\bf{Funding Statement}} This paper is supported by the National Natural Science Foundation of China (Nos.11901064 and 12071325).\\
{\bf{Conflict of interest}} All the authors in the paper have no conflict of
interest.\\
{\bf{Informed consent}} Informed consent was obtained from all individual
participants included in the study.
\section*{Author Contribution Statement} Dr. Peng He wrote the main manuscript text. Dr.
Xue-ping Wang rewrote and reviewed the whole article.

\end{document}